\documentclass[10pt]{amsart}
\usepackage{amssymb, amsthm, amsmath}
\usepackage{latexsym}

\theoremstyle{plain}
\newtheorem{prop}{Proposition}
\newtheorem{thm}{Theorem}
\newtheorem{cor}{Corollary}
\newtheorem{lemma}{Lemma}
\newtheorem*{movlemma}{Moving Lemma}

\theoremstyle{remark}

\theoremstyle{definition}
\newtheorem{defn}{Definition}

\def\AA{{\mathbb A}}
\def\bP{{\mathbb P}}

\def\Z{{\mathbb Z}}

\def\C{{\mathbb C}}
\def\cO{{\mathcal O}}
\def\F{{\mathcal F}}

\def\E{{\mathcal E}}
\def\cE{{\mathcal E}}
\def\cC{{\mathcal C}}
\def\cD{{\mathcal D}}

\def\D{{\mathcal D}}

\def\DS{{\mathfrak D}}
\def\QQ{{\mathcal Q}}
\def\ev{{\mathrm{ev}}}
\def\WW{{\widehat W}}
\def\l{{\lambda}}
\def\m{{\mu}}
\def\n{{\nu}}
\def\L{{\Lambda}}

\def\X{{\mathfrak X}}
\def\Y{{\mathfrak Y}}
\def\a{{\alpha}}

\def\s{{\sigma}}
\def\t{{\tau}}

\DeclareMathOperator{\Spec}{Spec}
\DeclareMathOperator{\rk}{rk}

\DeclareMathOperator{\Span}{Span}

\DeclareMathOperator{\Supp}{Supp}

\newcommand{\ssm}{\smallsetminus}
\newcommand{\gequ}{\geqslant}
\newcommand{\lequ}{\leqslant}
\newcommand{\lra}{\longrightarrow}

\newcommand{\ra}{\rightarrow}

\newcommand{\ov}{\overline}

\newcommand{\wt}{\widetilde}
\newcommand{\Pf}{\mbox{Pfaffian}}
\newcommand{\wh}{\widehat}

\begin{document}

\title[Quantum Cohomology of Orthogonal Grassmannians]
{Quantum Cohomology of Orthogonal Grassmannians}
\author{Andrew Kresch and Harry Tamvakis}
\date{March 10, 2003}
\address{Department of Mathematics, University of Pennsylvania,
209 South 33rd Street,
Philadelphia, PA 19104-6395, USA}
\email{kresch@math.upenn.edu}
\address{Department of Mathematics, Brandeis University, MS 050,
P.\ O.\ Box 9110, Waltham, MA 02454-9110, USA}
\email{harryt@brandeis.edu}

\begin{abstract}
Let $V$ be a vector space with a nondegenerate symmetric form
and $OG$ be the orthogonal
Grassmannian which  parametrizes maximal isotropic
subspaces in $V$. We give a presentation for the (small) quantum
cohomology ring $QH^*(OG)$ and show that
its product structure is determined by the ring of 
$\wt{P}$-polynomials. A 
`quantum Schubert calculus' is formulated,
which includes quantum Pieri
and Giambelli formulas, as well as algorithms for computing 
Gromov--Witten invariants. As an application, we show that the
table of $3$-point, genus zero Gromov--Witten invariants for $OG$ 
coincides with that for a corresponding Lagrangian Grassmannian $LG$, 
up to an involution.
\end{abstract}

\maketitle

\section{Introduction}
\label{intro}

\noindent
Consider a complex
vector space $V$ together with a nondegenerate symmetric form.
Our aim is to study the structure of the small quantum cohomology
ring of the orthogonal Grassmannian of maximal isotropic subspaces
in $V$. In a companion paper to this one \cite{KTlg}, we provide a 
similar analysis in type $C$, i.e., for the Lagrangian Grassmannian,
and the reader is referred there and to \cite{FP} \cite{LT}
for further background. The story
in the orthogonal case is similar, but with significant differences,
both in the results and in their proofs.

Assuming the dimension of $V$ is even and equals $2n+2$ for some 
natural number $n$, then the space of maximal isotropic
subspaces of $V$ has two connected components, each isomorphic to
the {\em even orthogonal Grassmannian} or {\em spinor variety}
$OG=OG(n+1,2n+2)=SO_{2n+2}/P_{n+1}$. Here $P_{n+1}$ is the 
maximal parabolic subroup of $SO_{2n+2}$
associated to a `right end root' in the 
Dynkin diagram of type $D_{n+1}$. 
We note that $OG(n+1,2n+2)$ is isomorphic (in fact
projectively equivalent) to the odd orthogonal
Grassmannian $OG(n,2n+1)=SO_{2n+1}/P_n$. Therefore it
suffices to work only with the even orthogonal example, 
and we will do so throughout this paper. We agree that 
a class $\a$ in the cohomology $H^{2k}(\X,\Z)$ of a complex variety 
$\X$ has degree $k$, to avoid doubling of all degrees.

The cohomology ring $H^*(OG,\Z)$ 
has a $\Z$-basis of Schubert classes $\t_{\l}$, 
one for each strict partition $\l=(\l_1>\l_2>\cdots>\l_{\ell}>0)$ with
$\l_1\lequ n$. Their multiplication can be described using the 
{\em $\wt{P}$-polynomials} of Pragacz and Ratajski
\cite{PR}. Let $X=(x_1,\ldots,x_n)$ be an $n$-tuple of
variables and define $\wt{P}_0(X)=1$ and 
$\wt{P}_i(X)=e_i(X)/2$ for each $i>0$, where
$e_i(X)$ denotes the
$i$-th elementary symmetric polynomial in $X$.
For nonnegative integers  $i,j$  with $i\gequ j$, set
\begin{equation}
\label{clgiam1}
\wt{P}_{i,j}(X)=
\wt{P}_i(X)\wt{P}_j(X)+
2\sum_{k=1}^{j-1}(-1)^k\wt{P}_{i+k}(X)\wt{P}_{j-k}(X)+(-1)^j
\wt{P}_{i+j}(X),
\end{equation}
and for any partition
$\l$ of length $\ell=\ell(\l)$, not necessarily strict, define
\begin{equation}
\label{clgiam2}
\wt{P}_{\l}(X)=\Pf[\wt{P}_{\l_i,\l_j}(X)]_{1\lequ i<j\lequ r},
\end{equation}
where $r=2\lfloor (\ell+1)/2\rfloor$. Let
$\D_n$ be the set of strict partitions $\l$ with $\l_1\lequ n$.

Let $\L'_n$ denote the $\Z$-algebra generated
by the polynomials $\wt{P}_{\l}(X)$ for all 
$\l\in\D_n$; $\L'_n$ is isomorphic
to the ring $\Z[X]^{S_n}$ of symmetric polynomials in $X$.
By results of \cite[Sect.\ 6]{P} and \cite{PR} we have that the map
sending $\wt{P}_{\l}(X)$
to $\t_{\l}$ for all $\l\in\D_n$ extends to a 
surjective {\em ring} homomorphism $\phi:\L_n'\ra H^*(OG,\Z)$ 
with kernel generated by the relations $\wt{P}_{i,i}(X)=0$ 
for $1\lequ i\lequ n$. The map $\phi$ can be realized as 
evaluation on the Chern roots of the tautological quotient 
vector bundle $Q$ over $OG$ (note that the top Chern class of $Q$
vanishes). In this way we obtain a presentation for 
the cohomology ring of $OG$, and 
equations (\ref{clgiam1}) and (\ref{clgiam2})
become Giambelli-type formulas, which express the Schubert classes
in terms of the special ones. 

We present an extension of these results to
the (small) quantum cohomology ring of $OG$, denoted $QH^*(OG)$. 
This is an algebra over $\Z[q]$, where $q$ is a formal
variable of degree $2n$ (the classical formulas are recovered
by setting $q=0$).

\begin{thm}
\label{ogthm} 
The map which sends
$\wt{P}_{\l}(X)$ to $\t_{\l}$ for all $\l\in\D_n$
and $\wt{P}_{n,n}(X)$ to $q$ extends to 
a surjective {\em ring} homomorphism $\L'_n\ra QH^*(OG)$
with kernel generated by the relations $\wt{P}_{i,i}(X)=0$ for 
$1\lequ i\lequ n-1$.
The ring $QH^*(OG)$ is presented as a quotient of the polynomial
ring $\Z[\t_1,\ldots,\t_n,q]$ modulo the relations 
\begin{equation}
\label{ogqrel}
\t_i^2+2\sum_{k=1}^{i-1}(-1)^k\t_{i+k}\t_{i-k}+(-1)^i\t_{2i}=0
\end{equation}
for all $i<n$, together with the quantum relation
\begin{equation}
\label{ogqreln}
\t_n^2=q
\end{equation}
{\em(}it is understood that $\t_j=0$ for $j>n${\em)}.
The Schubert class $\t_{\l}$
in this presentation is given by the Giambelli formulas
\begin{equation}
\label{ogqgiam1}
\t_{i,j}=
\t_i\t_j+2\sum_{k=1}^{j-1}(-1)^k\t_{i+k}\t_{j-k}+(-1)^j\t_{i+j}
\end{equation}
for $i>j>0$, and 
\begin{equation}
\label{ogqgiam2}
\t_{\l}=\text{\em Pfaffian}[\t_{\l_i,\l_j}]_{1\lequ i<j\lequ r},
\end{equation}
where quantum multiplication is employed
throughout. In other words, classical Giambelli and quantum Giambelli
coincide for $OG$.
\end{thm}
\noindent
We remark that the statements in Theorem \ref{ogthm}
are direct analogues of the corresponding facts for 
$SL_N$-Grassmannians \cite{Ber}. However, these results
stand in contrast to the case of the Lagrangian Grassmannian $LG(n,2n)$,
where quantum Giambelli does not coincide with classical
Giambelli on $LG(n,2n)$ (see \cite{KTlg} for more details).

Our proof of Theorem \ref{ogthm} follows the scheme of
\cite{KTlg}, with two main differences. We require a 
Pfaffian identity for type $D$ Schubert polynomials \cite[\S 3.3]{KTlg}, 
which gives a key relation in the Chow group of a certain
{\em orthogonal Quot scheme} $OQ_d$. The latter scheme 
compactifies the moduli space of degree $d$ maps $\bP^1\ra OG$;
however our definition of $OQ_d$ differs from that in the
Lagrangian case of \cite{KTlg}, as the direct analogue
of the Grothendieck Quot scheme \cite{G} here is 
not suitable for doing computations. 

In $QH^*(OG)$ there are formulas
\[
\t_{\l}\cdot \t_{\m} = 
\sum \langle \t_{\l}, \t_{\m}, \t_{\wh{\n}} \rangle_d \,\t_{\n}\, q^d,
\]
where the sum is over $d\gequ 0$ and strict partitions $\n$ with 
$|\n|=|\l|+|\m|-2nd$, and $\wh{\n}$ is the dual partition of $\n$,
whose parts complement the parts of $\n$ in the set
$\{1,\ldots,n\}$. Each quantum structure constant  
$\langle \t_{\l}, \t_{\m}, \t_{\wh{\n}} \rangle_d $
is a genus zero Gromov--Witten invariant for $OG$, and is a nonnegative 
integer. We present explicit 
formulas and algorithms to compute these numbers. This 
includes a quantum Pieri rule, which extends the classical
result of Hiller and Boe \cite{BH}. As an application, we show that
there is a direct identification between the $3$-point, genus zero
Gromov--Witten invariants
on $OG$ with corresponding ones for the Lagrangian Grassmannian
$LG(n-1,2n-2)$ (Theorem \ref{OGLG}).

This paper is organized as follows. In Section \ref{Ptilde} we study the 
$\wt{P}$-polynomials
and type $D$ Schubert polynomials, and prove a remarkable Pfaffian
identity for the latter. The orthogonal
Grassmannians are introduced in Section \ref{orthograss}, which includes 
a proof of the presentation for $QH^*(OG)$.
The proof of the quantum Giambelli formula (\ref{ogqgiam2})
of Theorem \ref{ogthm} is
done in Sections \ref{oquot} and  \ref{itoqd}, by studying intersections
on the orthogonal Quot scheme. In Section \ref{qsc} we
formulate a `quantum Schubert calculus' for $OG$. Finally,
the Appendix establishes an identity for $\wt{P}$-polynomials
which is used in \cite{KTdsp}.

The main results of this article and its companion paper \cite{KTlg} 
were announced at the Bonn Mathematische Arbeitstagung 2001 \cite{T}. 
The authors thank the Max-Planck-Institute f\"ur Mathematik
for its hospitality during the preparation of this paper. 
We also thank Anders Buch and Bill Fulton for useful correspondence.
Both authors were supported in part by National Science
Foundation post-doctoral research fellowships. 

\section{$\wt{P}$-polynomials and type $D$ Schubert polynomials}
\label{Ptilde}

\subsection{Basic definitions}
\label{bdefs}
All the notational conventions used in this section follow
\cite{KTdsp} and \cite{KTlg}. In particular, for strict
partitions $\l$ and $\m$, the difference $\l\ssm\m$ denotes the 
partition with parts given
by the parts of $\l$ which are not parts of $\m$.
A {\em composition} is a sequence of nonnegative integers with
only finitely many nonzero parts.
The $\wt{P}$-polynomials make sense when indexed
by any composition $\nu$, and satisfy Pfaffian relations
\begin{equation}
\label{pfaff}
\wt{P}_{\nu}(X)=
\sum_{j=1}^{g-1}(-1)^{j-1}\wt{P}_{\nu_j,\nu_g}(X)
\cdot \wt{P}_{\nu\ssm\{\nu_j,\nu_g\}}(X),
\end{equation}
where $g$ is an even number such that $\nu_i=0$ for $i>g$.
Define also the $\wt{Q}$-polynomial 
$\wt{Q}_{\n}(X)=2^{\ell}\,\wt{P}_{\n}(X)$ for each composition
$\n$ with $\ell$ nonzero parts. The $\wt{Q}$-polynomials have integer
coefficients, and span 
the ring $\Z[X]^{S_n}$ of symmetric functions in $n$
variables.

Let $\wt{W}_n$ be 
the Weyl group for the root system $D_n$, whose elements are denoted
as barred permutations. Recall that
$W_n$ is generated by the elements $s_{\Box},s_1,\ldots,s_{n-1}$: for $i>0$,
$s_i$ is the transposition interchanging $i$ and $i+1$, and 
$s_{\Box}$ is defined by
\[
(u_1,u_2,u_3,\ldots,u_n)s_{\Box}=(\ov{u}_2,\ov{u}_1,u_3,\ldots,u_n).
\]
Let $\wt{w}_0$ denote the element of maximal length in $\wt{W}_n$.
For each $\l\in \D_{n-1}$ we have a 
{\em maximal Grassmannian element} $w_{\l}$ of $\wt{W}_n$, defined as
in \cite[\S 3.2]{KTdsp}.

Each generator $s_i$ acts naturally on the polynomial ring $A[X]$,
where $A=\Z[\frac{1}{2}]$; for
$i>0$, $s_i$ interchanges $x_i$ and $x_{i+1}$, while $s_{\Box}$ 
sends $(x_1,x_2)$ to $(-x_2,-x_1)$; all other variables remain
fixed. There are divided difference operators $\partial'_i$ and 
$\partial_{\Box}$ on $A[X]$; for $i>0$ they are defined by
\[
\partial'_i(f)=(f-s_if)/(x_{i+1}-x_i)
\]
while 
\[
\partial_{\Box}(f)=(f-s_{\Box}f)/(x_1+x_2),
\]
for all $f\in A[X]$. These give rise to operators 
$\partial'_w: A[X]\ra A[X]$ for each element $w\in \wt{W}_n$, as
in \cite[\S 3.2]{KTdsp}. 

For all $w\in \wt{W}_n$ we have a {\em type $D$
Schubert polynomial} $\DS_w(X)\in A[X]$ defined by
\[
\DS_w(X)=(-1)^{n(n-1)/2}\partial'_{w^{-1}\wt{w}_0}
\Bigr(x_1^{n-1}x_2^{n-2}\cdots x_{n-1}\wt{P}_{n-1}(X)\Bigl).
\]
These type $D$ polynomials were defined in \cite[\S 3.3]{KTdsp}; they
agree with the orthogonal Schubert
polynomials of \cite{LP2} up to a sign, which depends on the degree.
The polynomial $\DS_w(X)$ 
represents the Schubert class associated to $w$ in the cohomology
ring of the flag manifold $SO_{2n}/B$. Let us define 
$\DS_{\l}'(X)=\DS_{w_{\l}s_{\Box}}(X)$. It follows from the definitions
and \cite[Theorem 7]{KTdsp} that 
$\DS_{\l}'(X)=\partial_{\Box}(\wt{P}_{\l}(X))$,
for all non-zero partitions $\l\in\D_{n-1}$.

\subsection{A Pfaffian identity}
\label{pfids}
We require the identity in the following theorem for our
proof of the quantum Giambelli formula for $OG(n+1,2n+2)$. 

\begin{thm}
\label{pfid1} 
Fix $\l\in\D_n$ of length $\ell\gequ 3$, and set
$r=2\lfloor(\ell+1)/2\rfloor$. Then 
\begin{equation}
\label{pf1}
\sum_{j=1}^{r-1}(-1)^{j-1}\,\DS'_{\l_j,\l_r}(X)\,
\DS'_{\l\ssm\{\l_j,\l_r\}}(X)=0.
\end{equation}
\end{thm}

\begin{proof}
We first observe, using the homogeneity of the two sides,
that (\ref{pf1}) is equivalent to the identity 
\begin{equation}
\label{Qtildeid}
\sum_{j=1}^{r-1}(-1)^{j-1}\,\partial_{\Box}(\wt{Q}_{\l_j,\l_r}(X))\cdot
\partial_{\Box}(\wt{Q}_{\l\ssm\{\l_j,\l_r\}}(X))=0
\end{equation}
for $\wt{Q}$-polynomials, 
which should hold for $\l$ and $r$ as in the theorem.

Let $X''=(x_3,\ldots,x_n)$ and define
\[
m_{r,s}(x_1,x_2)=
\begin{cases}
x_1^rx_2^s+x_1^sx_2^r&\text{if $r\neq s$},\\
x_1^rx_2^r&\text{if $r=s$}
\end{cases}
\]
to be the monomial symmetric function in $x_1$ and $x_2$. For any 
partition $\l$ and nonnegative integers $a$ and $b$, let $C(\l,a,b)$ 
denote
the set of compositions $\mu$ with $\l_i-\mu_i\in\{0,1,2\}$ for 
all $i$ and $\l_i-\m_i=1$ (resp.\ $\l_i-\m_i=2$)
for exactly $a$ (resp.\ $b$) values of $i$.
\begin{prop}
\label{Boxprop}
For any nonzero strict partition $\l$, we have
\begin{equation}
\label{Bact}
\partial_{\Box}(\wt{Q}_{\l}(X)) =2
\sum_{\substack{0\lequ s\lequ r \lequ \ell \\ r+s \ \mathrm{even}}} 
m_{r,s}(x_1,x_2) 
\sum_{\substack{a+2b=r+s+1 \\ 0\lequ b\lequ s}}
\binom{a-1}{s-b}\sum_{\mu\in C(\l,a,b)} \wt{Q}_{\mu}(X'').
\end{equation}
\end{prop}
\begin{proof} Let $X'=(x_2,\ldots,x_n)$. According to
\cite[Prop.\ 1]{KTlg}, for any partition $\l$ of length $\ell$
(not necessarily strict), we have
\begin{equation}
\label{extnone}
\wt{Q}_{\l}(X)=
\sum_{k=0}^{\ell}x_1^k\sum_{\mu\in B(\l,k)}\wt{Q}_{\mu}(X'),
\end{equation}
where $B(\l,k)$ is defined to be the set of all compositions $\mu$ such that 
$|\l|-|\mu|=k$ and $\l_i-\mu_i\in\{0,1\}$ for each $i$. By 
applying (\ref{extnone}) twice we obtain
\begin{equation}
\label{extntwo}
\wt{Q}_{\l}(X)=
\sum_{0\lequ s \lequ r \lequ \ell}m_{r,s}(x_1,x_2)
\sum_{\substack{j+2k=r+s\\ 0\lequ k\lequ s}}
\binom{j}{s-k}\sum_{\mu\in C(\l,j,k)}\wt{Q}_{\mu}(X'').
\end{equation}

Suppose that $r\gequ s\gequ 0$. If $r+s$ is even, then
$\partial_{\Box}(m_{r,s}(x_1,x_2))=0$. If $r+s$ is odd, we have
\[
\partial_{\Box}(m_{r,s}(x_1,x_2))=
2\sum_{\substack{c+d=r+s-1\\ c,d\gequ s}}(-1)^{c-s}x_1^cx_2^d.
\]
We now apply this to (\ref{extntwo}) and gather terms to obtain (\ref{Bact}).
\end{proof}

\medskip

\noindent
{\bf Example.} For  all $a$, $b$ with  $a>b\gequ 0$, we have
\begin{align}
\label{exequ}
\begin{split}
\partial_{\Box} ( \wt{Q}_{a,b} (X)) &=
2 \,\Bigr(\wt{Q}_{a-1,b}(X'')+\wt{Q}_{a,b-1}(X'')\Bigl) \medskip \\
 &+2\,x_1x_2\Bigr( \wt{Q}_{a-2,b-1}(X'')+\wt{Q}_{a-1,b-2}(X'') \Bigl).
\end{split}
\end{align}
In the equation (\ref{exequ}) and later on we agree that
$\wt{Q}_{\mu}(X'')=0$ if any of the components of $\mu$ 
are negative.

\medskip

As in \cite[\S 2.3]{KTlg}, the rest of the argument can be expressed
using only the partitions which index the polynomials involved. We
thus begin by defining a commutative
$\Z$-algebra  ${\mathcal B}$ with formal variables which 
represent these indices. The algebra ${\mathcal B}$ is
generated by symbols $(a_1,a_2,\ldots)$, where the entries $a_i$ are 
barred integers; each $a_i$ can have up to two bars.
The symbol $(a_1,a_2,\ldots)$ corresponds to
the polynomial $\wt{Q}_{\m}(X'')$, where $\m$ is the composition
with $\m_i$ equal to the integer $a_i$ minus the number of bars
over $a_i$. We identify $(a,0)$ with $(a)$.

Let $\mu$ be a barred partition, that is, a partition in which bars
have been added to some of the entries. For $\ell(\m)\gequ 3$,
we impose the Pfaffian relation
\begin{equation}
\label{pfrl}
(\m)=
\sum_{j=1}^{m-1}(-1)^{j-1}(\m_j,\m_m)
\cdot (\m\ssm\{\m_j,\m_m\}),
\end{equation}
which corresponds to (\ref{pfaff}) for $\n=\m$ (here $m=2\lfloor 
(\ell(\m)+1)/2\rfloor$, as usual).
Iterating this gives
\begin{equation}
\label{defin}
(\m)=\sum \epsilon(\m,\nu)(\nu_1,\nu_2)\cdots(\nu_{m-1},\nu_m),
\end{equation}
where the sum is over all $(m-1)(m-3)\cdots (1)$ ways to write the
set $\{\m_1,\ldots,\m_m\}$ as a union of pairs
$\{\nu_1,\nu_2\}\cup\cdots\cup\{\nu_{m-1},\nu_m\}$, and where 
$\epsilon(\m,\nu)$ is the sign of the permutation that takes
$(\m_1,\ldots,\m_m)$ into $(\nu_1,\ldots,\nu_m)$; we adopt the
convention that $\nu_{2i-1}\gequ\nu_{2i}$.

We also define the square bracket symbols $[a]=(\ov{a})$ and 
$[a,b]=(\ov{a},b)+(a,\ov{b})$,
where $a$ and $b$ are integers, each with up to one bar. For example,
the right hand side of equation (\ref{exequ}) corresponds to the
sum $2\,[a,b]+2\,x_1x_2\,[\ov{a},\ov{b}]$ in ${\mathcal B}[x_1,x_2]$.
Finally, we impose the relations
\begin{equation}
\label{check}
[a,b]=(\ov{a})(b)-(a)(\ov{b}) 
\end{equation}
for integers $a$, $b$, with up to one bar each; this agrees with
a corresponding identity 
\[
\wt{Q}_{a-1,b}+\wt{Q}_{a,b-1}=\wt{Q}_{a-1}\wt{Q}_b-
\wt{Q}_a\wt{Q}_{b-1}
\]
of $\wt{Q}$-polynomials.

Using these conventions and equations (\ref{Bact}) and (\ref{exequ}), 
we are reduced to showing that $S_1+S_2=0$, where
\[
S_1=\sum_{\substack{a+2b=r+s+1 \\ 0\lequ b\lequ s}}\binom{a-1}{s-b}
\sum_{j=1}^{r-1}(-1)^{j-1}\,[\l_j,\l_r]
\sum_{\mu\in C(\l\ssm\{\l_j,\l_r\},a,b)} (\m),
\]
\[
S_2=\sum_{\substack{a'+2b'=r+s-1 \\ 0\lequ b'\lequ s-1}}\binom{a'-1}{s-b'-1}
\sum_{j=1}^{r-1}(-1)^{j-1}\,[\ov{\l}_j,\ov{\l}_r]
\sum_{\mu\in C(\l\ssm\{\l_j,\l_r\},a',b')} (\m),
\]
and  $r\gequ s\gequ 0$ are fixed integers with $r+s$ even. The proof
of this is rather similar to the proofs of Theorems 2 and 3 of \cite{KTlg},
and we will point out only the main difference here.

We first apply (\ref{defin}) to expand the terms $(\mu)$ in both $S_1$
and $S_2$. The cancellation technique of 
\cite[\S 2.3]{KTlg}, notably, the identity
\begin{equation}
\label{pf}
[a,b][c,d]-[a,c][b,d]+[a,d][b,c]=0,
\end{equation}
implies the vanishing of the sum of those
summands in $S_1$ which contain a pair with exactly
one bar, or at least two pairs with exactly three bars. The remainder
is a sum $S_1'$ consisting of those summands in $S_1$ with a unique
pair containing three bars, and no pair with only one bar. In the 
same way, one checks the vanishing of the sum of those
summands in $S_2$ which contain a pair with exactly
three bars, or at least two pairs with exactly one bar. There remains
a sum $S_2'$ consisting of those summands in $S_2$ with a unique
pair containing only one bar, and no pair with exactly three bars. 
Hence, it is enough to show that $S_1'+S_2'=0$. 

There is an obvious bijection between the summands in $S_1'$ and
$S_2'$, obtained by adding two bars to the unbarred part of the 
pair in $S_2'$ which contains only one bar (note that the corresponding
binomial
coefficients agree, as $(a,b)=(a',b'+1)$ for these two summands). To prove
that the sum of all corresponding terms is zero, it suffices to 
show that the expression
\begin{equation}
\label{needid}
\Bigl([a,b][\ov{c},\ov{d}]-[a,c][\ov{b},\ov{d}]+[a,d][\ov{b},\ov{c}]\Bigr)+
\Bigl([\ov{a},\ov{b}][c,d]-[\ov{a},\ov{c}][b,d]+[\ov{a},\ov{d}][b,c]\Bigr)
\end{equation}
vanishes identically in ${\mathcal B}$ (we then apply this with
$a=\l_r$, always). To check this, begin from the basic identities
\begin{equation}
\label{id1}
[a,b][\ov{c},\ov{d}]-[a,\ov{c}][b,\ov{d}]+[a,\ov{d}][b,\ov{c}]=0
\end{equation}
and
\begin{equation}
\label{id2}
[\ov{a},\ov{b}][c,d]-[\ov{a},c][\ov{b},d]+[\ov{a},d][\ov{b},c]=0
\end{equation}
which are easily shown using (\ref{check}).
Let $\left<x,y\right>=[\ov{x},y]+[x,\ov{y}]$ and note that
\begin{equation}
\label{id3}
\left<a,b\right>\left<c,d\right>-\left<a,c\right>
\left<b,d\right>+\left<a,d\right>\left<b,c\right>=0,
\end{equation}
which is shown using 
$\left<x,y\right>=(\ov{\ov{a}})(b)-(a)(\ov{\ov{b}})$
(another consequence of (\ref{check})).
The vanishing of (\ref{needid}) follows by combining (\ref{id1}),
(\ref{id2}) and (\ref{id3}).
\end{proof}

\section{Orthogonal Grassmannians}
\label{orthograss}

\subsection{Schubert varieties and incidence loci}
\label{schuvar}
Let $V$ be a fixed $(2n+2)$-dimensional complex vector space
equipped with a nondegenerate symmetric bilinear form on $V$.
The principal object of study is the orthogonal Grassmannian
$OG(n+1,2n+2)$ which is one component of the parameter space
of $(n+1)$-dimensional isotropic subspaces of $V$.
When $n$ is fixed, we write $OG$ for $OG(n+1,2n+2)$.
We have $\dim_{\C}OG=n(n+1)/2$. The identities in cohomology that we 
establish in this section remain valid if we work
over an arbitrary base field, and use Chow rings in place of cohomology.

Let $F_{\bullet}$ be a fixed complete isotropic flag of subspaces of $V$.
By convention, then, $OG$ parametrizes maximal isotropic spaces
$\Sigma\subset V$ such that $\Sigma\cap F_{n+1}$ has even codimension
in $F_{n+1}$.
We define the alternative flag $\wt{F}_\bullet$ to be the flag
$F_1\subset\cdots\subset F_n\subset \wt{F}_{n+1}$, where $\wt{F}_{n+1}$
is the unique maximal isotropic space containing $F_n$ but
not equal to $F_{n+1}$. We let
\begin{equation}
\label{defnfi}
F_\bullet^{(i)}=\begin{cases}
F_\bullet & \text{if $i \equiv (n+1)$ mod $2$},\\
\wt{F}_\bullet & \text{otherwise}.
\end{cases}
\end{equation}

The Schubert varieties $\X_\l\subset OG$
are indexed by partitions $\l\in \D_n$. We record
two ways to
write the conditions which define the Schubert variety $\X_\l$:
\begin{gather}
\X_\l = \{\,\Sigma\in OG\,|\,
\rk(\Sigma\to V/F_{n+1-\l_i})\lequ n+1-i,\,\,i=1,\ldots,\ell(\l)\,\}
\label{schubalt} \\
= \{\,\Sigma\in OG\,|\,
\rk(\Sigma\to V/F^{(i)\perp}_{n+1-\l_i})
\lequ n+1-i-\l_i,\,\,i=1,\ldots,\ell(\l)+1\,\}.
\label{schubert}
\end{gather}

Let $\t_{\l}$ be the class of $\X_{\l}$ in $H^*(OG,\Z)$.
The classical Giambelli formula (\ref{ogqgiam2})
for $OG$ is equivalent to the following identity in $H^*(OG,\Z)$:
\begin{equation}
\label{classicalgiambelli}
\t_\l = \sum_{j=1}^{r-1} (-1)^{j-1} \t_{\l_j,\l_r}\cdot
\t_{\l\smallsetminus\{\l_j,\l_r\}},
\end{equation}
for $r=2\lfloor (\ell(\l)+1)/2 \rfloor$. Let $\rho_n=(n,n-1,\ldots,1)$
and for $\m\in \D_n$, denote by $\wh{\mu}=\rho_n\ssm\mu$, the dual 
partition.
The Poincar\'e duality pairing on $OG$ satisfies
\[
\int_{OG} \t_{\l}\,\t_{\m} = \delta_{\l\wh{\m}}.
\]

Given an isotropic space $A\subset V$ of dimension $n-k$ ($k\gequ 0$),
the variety of
maximal isotropic spaces containing $A$ is a translate of the Schubert variety
$\X_{n,n-1,\ldots,k+1}$.
We have the following result on intersections of such varieties with the
Schubert varieties $\X_{\l}$; this is 
analogous to a similar result in type $C$ (\cite[Prop.\ 3]{KTlg}).

\begin{prop}
\label{meetsmallergrass}
Let $k\gequ 0$ and $\l\in \D_n$.
Let $A$ be an isotropic subspace of $V$ of dimension $n-k$, and 
let $Y\subset OG$ be the subvariety of
maximal isotropic subspaces of $V$ which contain $A$.
Then $\X_\lambda\cap Y$ is a Schubert variety in $Y\simeq OG(k+1,2k+2)$.
Moreover, if $\ell(\l)<k$ then the intersection,
if nonempty, has positive dimension.
\end{prop}

\begin{proof}
As in \cite{KTlg}, the intersection is defined by the attitude of
$\Sigma/A$ with respect to $F'_\bullet$, where
$F'_i=((F_i+A)\cap A^\perp)/A$.
For the intersection to be a point would require at least $k$
rank conditions, and hence $\ell(\l)\gequ k$.
\end{proof}

\medskip
The space $OG(n-1,2n+2)$ is the parameter space of lines on $OG$. For
a nonempty partition $\l$,
the variety of lines incident to $\X_\l$ 
is the Schubert variety $\Y_\l$, consisting of those
$\Sigma'\in OG(n-1,2n+2)$ such that
\begin{equation}
\label{schubertprime}
\rk(\Sigma'\to V/F^{(i)\perp}_{n+1-\l_i})\lequ n+1-i-\l_i,
\,\ \ \mathrm{for} \ \ i=1,\ldots,\ell+1.
\end{equation}
The codimension of $\Y_\l$ is $|\l|-1$.
Note that (i) the rank conditions (\ref{schubertprime})
are identical to those in (\ref{schubert});
(ii) the rank condition corresponding to $i=\ell(\l)+1$, which was 
redundant in defining the Schubert varieties in $OG$, is necessary here.

\subsection{A Pfaffian identity on $OG(n-1,2n+2)$}
\label{pfaffident}
Let $F=F_{SO}(V)$ denote the variety of complete isotropic flags
in $V=\C^{2n+2}$. There is a natural projection map from $F$ to
the orthogonal Grassmannian $OG(n-1,2n+2)$, inducing an
injective pullback morphism on cohomology. Introduce an extra
variable $x_{n+1}$ and let $X^+=(x_1,\ldots,x_{n+1})$.
Referring to 
\cite[\S 2.4 and Sect.\ 3]{KTdsp}, one checks that the Schubert class
$[\Y_\l]$ in $H^*(OG(n-1,2n+2))$ pulls back to the class represented
by $\DS'_{\l}(X^+)$ in $H^*(F)$, for each $\l\in\D_{n-1}$. 
Here $X^+$ corresponds to the vector of Chern roots of the dual
to the tautological rank $n+1$ vector bundle over $F$, ordered as
in \cite[Sect.\ 2]{KTdsp}. Theorem \ref{pfid1} remains true with
$X^+$ in place of $X$, and gives

\begin{cor}
\label{forcor} 
For every $\l\in\D_n$ of length $\ell\gequ 3$ 
and $r=2\lfloor(\ell+1)/2\rfloor$ we have
\begin{equation}
\label{idone}
\sum_{j=1}^{r-1}(-1)^{j-1}\,[\Y_{\l_j,\l_r}]\,
[\Y_{\l\ssm\{\l_j,\l_r\}}]=0
\end{equation}
in $H^*(OG(n-1,2n+2),\Z)$.
\end{cor}

\subsection{Quantum relations and two-condition Giambelli}

Recall that in $QH(OG)$, the degree of $q$ is
$$\int_{OG} c_1(T_{OG})\cdot \t_{\wh{1}}=2n.$$
It follows, for degree reasons, that
the relations in cohomology (\ref{ogqrel}) and
the quantum Giambelli formula for the two-condition Schubert
classes (\ref{ogqgiam1}) -- which
we know to hold classically -- hold in $QH(OG)$.
The degree $2n$ quantum relation (\ref{ogqreln}) follows from the
elementary enumerative fact that there is a unique line on $OG$
through a given point, incident to two general translates of $\X_n$.
Arguing as in \cite{ST}, now, we obtain a presentation of $QH^*(OG)$ 
as a quotient of the polynomial ring $\Z[\t_1,\ldots,\t_n,q]$ modulo the
relations (\ref{ogqrel}) and (\ref{ogqreln}) (see also \cite[Sect.\ 10]{FP}).

The proof of the more difficult quantum Giambelli formula (\ref{ogqgiam2})
occupies Sections \ref{oquot} and \ref{itoqd}.

\section{Orthogonal Quot schemes}
\label{oquot}

\subsection{Overview}
\label{overview}

In the next two sections,
we define the orthogonal Quot scheme and establish an identity in
its Chow group,
from which identity (\ref{ogqgiam2}) in $QH^*(OG)$ readily follows.
We make use of type $D$ degeneracy loci for isotropic morphisms
of vector bundles \cite{KTdsp} to define classes
$[W_\l(p)]_k$ ($p\in \bP^1$)
of the appropriate dimension $k:=n(n+1)/2+2nd-|\l|$ in the
Chow group of the
orthogonal Quot scheme $OQ_d$, which compactifies the space of 
degree-$d$ maps
$\bP^1\to OG$.
Let $p'\in\bP^1$ be distinct from $p$, and denote by $W'$
the degeneracy locus defined by a general translate of the fixed
isotropic flag $F_\bullet$.
We produce a Pfaffian formula analogous to (\ref{classicalgiambelli}):
\begin{equation}
\label{Wequ}
[W_{\l}(p)]_k=\sum_{j=1}^{r-1}(-1)^{j-1}[W_{\l_j,\l_r}(p)\cap
W'_{\l\ssm\{\l_j,\l_r\}}(p')]_k,
\end{equation}
for any $\l\in\D_n$ with $\ell(\l)\gequ 3$ and
$r=2\lfloor(\ell(\l)+1)/2\rfloor$.

As in \cite{KTlg}, we need the cycles in (\ref{Wequ})
to remain rationally equivalent under
further intersection with some (general translate of) $W_{\m}(p'')$,
for $\m\in \D_n$ and $p''\in \bP^1$ distinct from $p$, $p'$,
Also, as in loc.\ cit., we accomplish this by working on
a modification $OQ_d(p'')$, on which the
evaluation-at-$p''$ map is globally defined, and employing
refined intersection operation from $OG$.

The rational equivalences that we produce --- (\ref{Wequ}) and a similar
equivalence on $OQ_d(p'')$ --- come by combining equivalences of the
following types:
(i) the classical Pfaffian formulas on $OG$ (\ref{classicalgiambelli});
(ii) the Pfaffian identities (\ref{idone}) on $OG(n-1,2n+2)$;
(iii) rational equivalences $\{p\}\sim\{p'\}$ on $\bP^1$.
Indeed, the essence of (iii) is that we can replace
$p'$ with $p$ in (\ref{Wequ});
the intersection
$W_{\l_j,\l_r}(p)\cap
W'_{\l\ssm\{\l_j,\l_r\}}(p)$
now has $k$-dimension components supported in the boundary of the
Quot scheme.
The cancellation of these contributions in the Chow group is
precisely equation (\ref{idone}).

\subsection{Definition of $OQ_d$}
\label{lagrquotdefn}

Let $V$ be a complex vector space $V$ of dimension $N=r+s$ and
fix $d\gequ 0$. Following Grothendieck \cite{G}, 
there is a smooth projective variety $Q_d$, the
{\em Quot scheme}, which parametrizes flat families of quotient
sheaves of $\cO_{\bP^1}\otimes V$ with Hilbert polynomial $p(t)=st+s+d$.
This variety
compactifies the space of parametrized degree-$d$ maps from $\bP^1$
to the Grassmannian of $r$-dimensional subspaces of $V$.
On $\bP^1\times Q_d$ there is a universal exact sequence of sheaves
\begin{equation}
\label{sheafseq}
0\lra \E\lra \cO\otimes V\lra \QQ\lra 0
\end{equation}
with $\E$ locally free of rank $r$.  From now on, we fix $V$ as
in Section \ref{orthograss} and $r=s=n+1$.

\begin{defn}
Let $d$ be a nonnegative integer.
The {\em isotropic locus} $Q_d^{\rm iso}$ is
the closed subscheme of $Q_d$ which is defined by the vanishing
of the composite
$$\E\lra \cO_{\bP^1}\otimes V 
\stackrel{\alpha}\lra \cO_{\bP^1}\otimes V^*\lra \E^*$$
where $\alpha$ is the isomorphism defined by the given bilinear form on $V$.
\end{defn}

The embedding of $OG$ in the Grassmannian $G(n+1,2n+2)$ of
$(n+1)$-dimensional subspaces of $V$ is degree-doubling, that is,
in the sheaf sequence (\ref{sheafseq}) corresponding to degree-$d$ maps
$\bP^1\to OG$, the sheaf $\QQ$ has degree $2d$.
For any $d$, $Q_{2d}^{\rm iso}$ contains an open subscheme isomorphic
to the moduli space $M_{0,3}(OG,d)$:
\begin{defn}
Let $d$ be a nonnegative integer.
Then $OM_d$ is the open subscheme of $Q_{2d}^{\rm iso}$ defined by
the conditions (i) $\E\to \cO_{\bP^1}\otimes V$ has everywhere full rank;
(ii) the image of $\E\to \cO_{\bP^1}\otimes V$ at any point has
intersection with $F_{n+1}$ of dimension congruent to $(n+1)$ mod $2$.
\end{defn}

Unfortunately, $Q_{2d}^{\rm iso}$ generally has components of dimension larger
than the dimension of $OM_d$.
The remedy is to throw away any point of (\ref{sheafseq}) where the
rank of $\cE\to \cO\otimes V$ drops by just $1$ at some point of
$\bP^1$.
We can do this, and still be left with a closed subscheme of
$Q_{2d}^{\rm iso}$, because in any degeneration situation in which
the rank of $\cE\to \cO\otimes V$ drops from full to less than full,
the drop is by at least $2$.

\begin{defn}
For $d\in (1/2)\Z$,
the {\em orthogonal Quot scheme} $OQ_d$ is the subset of
$Q_{2d}^{\rm iso}$ consisting of points whose sheaf sequence
(\ref{sheafseq}) satisfies
$\rk(\cE_p\to V)\ne n$ for all $p\in \bP^1$,
and such that where it has full rank, the image has intersection
with $F_{n+1}$ of even codimension in $F_{n+1}$.
This subset, evidently constructible and closed by virtue of
Proposition \ref{closedunderspec}, below,
is given the reduced scheme structure.
\end{defn}

\begin{lemma}
\label{grassfiber}
Let $\psi\colon C_0\to G(n+1,2n+2)$ be a morphism, with $C_0\cong \bP^1$,
and let $C$ be a tree of $\bP^1$'s containing $C_0$ and
$\varphi\colon C\to G(n+1,2n+2)$ a map which restricts to $\psi$ on $C_0$.
Let $$\widetilde C:=C_1\cup C_2\cup\cdots \cup C_m$$
{\em(}$m\gequ 1${\em)} denote a chain of components in $C$,
with $C_i\ne C_0$ for all $i\gequ 1$,
and assume $C_1$ meets $C_0$ at the point $p$ and $C_i$ is collapsed
by $\varphi$ for all $i$ with $1\lequ i\lequ m-1$.
Let $\pi\colon C\to C_0$ denote the morphism which collapses all components
of $C$ except $C_0$.
Let
$$0\to \cE_0\to \cO\otimes V\to \QQ_0\to 0$$
denote the pullback of the universal sequence via $\psi$,
and let
$$0\to \cE\to \cO\otimes V\to \QQ\to 0$$
denote the pullback of the universal sequence via $\varphi$
{\em(}so that $\cE|_{C_0}\simeq \cE_0${\em)}.
Assume the restriction of $\cE$ to $C_m$ splits as
$$\cO(-b_1)\oplus \cdots\oplus \cO(-b_j)\oplus \cO^{n+1-j}$$
with $b_1, \ldots, b_j\gequ 1$.
Then the morphism $\pi_*\cE\to \pi_*(\cO\otimes V)=\cO\otimes V$
factors through $\cE_0$,
and the cokernel of $\pi_*\cE\to\cE_0$ is a torsion sheaf whose
fiber at $p$ has dimension at least $j$.
\end{lemma}

\begin{proof}
We may choose $n-j$ independent sections $s_1$, $\ldots$, $s_{n-j}$ of
$\cE|_{C_m}$.
These extend uniquely to $n-j$ independent sections of
$\cE|_{\widetilde C}$, and hence span an
$(n-j)$-dimensional subspace $\Sigma$ of the fiber of $\cE$ at the point $p$.
The map $(\pi_*\cE)_p\to(\cE_0)_p$ on fibers at $p$ has image contained
in $\Sigma$. Hence the dimension of the fiber at $p$ of the 
cokernel of $\pi_*\cE\to\cE_0$ is at least $j$.
\end{proof}

\begin{prop}
\label{closedunderspec}
For any $d\in (1/2)\Z$, the subset $OQ_d\subset Q_{2d}^{\rm iso}$ is
closed under specialization.
\end{prop}

\begin{proof}
Suppose $x_1\in OQ_d$ specializes to $x_0\in Q_{2d}$.
Then there is a discrete valuation ring $R$ and a morphism
$\varphi\colon \Spec R\to Q_{2d}$ such that
the generic point maps to $x_1$ and the special point maps to $x_0$.

Denote the fraction field of $R$ by $K$ and the residue field by $k$.
It suffices to consider the case where $x_0$ is a closed point,
hence $k=\C$ is algebraically closed.
We show that given
the exact sequence of coherent sheaves at the generic point
\begin{equation}
\label{seqgeneric}
0\to \cE\to \cO\otimes V\to \QQ\to 0
\end{equation}
on $\bP^1_K$, we can reconstruct the map $\varphi$ and hence the
sheaf sequence at the special point (possibly
replacing $R$ by its integral closure in a finite extension
of $K$).
Then, we note that the torsion of the quotient sheaf at the special point
cannot have rank $1$ at any point of $\bP^1_k$.

Let the sequence (\ref{seqgeneric}) be given.
The support of $\QQ^{\rm tors}$ specializes to a well-defined closed
subset $Z\subset\bP^1_k$; we let $Y=\Supp(\QQ^{\rm tors})\cup Z$.
Now consider:
\begin{equation}
\label{seqtfree}
0\to \cE'\to \cO\otimes V\to \QQ/\QQ^{\rm tors}\to 0
\end{equation}
on $\bP^1_K$.
This corresponds to a morphism $\bP^1_K\to OG$ (the actual map to the
orthogonal Grassmannian
underlying the sheaf sequence (\ref{seqgeneric})).
By replacing $K$ by a finite extension and $R$ by its integral closure
in the extension, if necessary, then there exists, by semistable reduction,
a modification
$$\pi\colon S\to \bP^1_R$$
with exceptional divisor a tree of $\bP^1$'s, and a
morphism
$S\to OG$,
such that $\pi$ restricts to the given morphism
$\bP^1_K\to OG$.  We consider the pullback of the
universal exact sequence
$$0\to\widetilde \cE\to \cO\otimes V\to \widetilde \QQ\to 0$$
on $S$.
Pushing forward the map $\cE\to \cO\otimes V$ by $\pi$ yields an exact sequence
\begin{equation}
\label{seqspecial}
0\to \pi_*\widetilde\cE\to \cO\otimes V\to \cC\to 0
\end{equation}
The cokernel $\cC$, being a subsheaf of $\pi_*\widetilde\QQ$,
is torsion-free over $\Spec R$, and hence flat:
(\ref{seqspecial}) corresponds to the map from $\Spec R$ to the (possibly
smaller degree) Quot scheme determined by (\ref{seqtfree}).

We extend (\ref{seqgeneric}) to all of $\bP^1_R$
by patching and pushing forward.
The sequences (\ref{seqgeneric}) on $\bP^1_K$ and
(\ref{seqspecial}) on $\bP^1_R\smallsetminus Y$ patch to give the
sequence
$$0\to\widehat\cE\to\cO\otimes V\to\widehat\QQ\to 0$$
on $\bP^1_R\smallsetminus Z$.
Pushing forward via $i\colon \bP^1_R\smallsetminus Z\to \bP^1_R$ gives
\begin{equation}
\label{specializesto}
0\to i_*\widehat\cE\to \cO\otimes V\to\cD\to 0,
\end{equation}
(where $\cD$ is the indicated cokernel), flat 
over $\bP^1_R$ since $i_*\widehat\cE$ is locally free.
This gives the morphism $\varphi\colon \Spec R\to Q_{2d}$
that we started with.

We now consider the restriction of (\ref{specializesto}) to the special fiber:
$$0\to (i_*\widehat\cE)_k\to \cO\otimes V\to \cD_k\to 0,$$
and verify it satisfies the rank conditions.
By semicontinuity, the dimension of the fiber of
$\cD_k^{\rm tors}$ is $\gequ 2$ at every
point of $Z$.
Suppose $p$ is a point in $\bP^1_k\smallsetminus Z$.
Then $\cD_k$, on a neighborhood of $p$, is isomorphic to
$\cC_k:=\cC\otimes_Rk$,
so it suffices to show every nonzero fiber of $\cC_k^{\rm tors}$ has
dimension $\gequ 2$.
Letting $(\,)_k$ denote restriction to the special fiber, we have:
$(\pi_*\widetilde\cE)_k\to \cO\otimes V$ factors through
$(\pi_k)_*(\widetilde\cE_k)\to\cO\otimes V$,
which in turn factors through a vector subbundle
$[(\pi_k)_*(\widetilde\cE_k)]'$
of $\cO\otimes V$
(the pullback of the universal subbundle by the actual map
$\bP^1_k\to OG$
at the special fiber),
and $\dim \cC_k^{\rm tors}\otimes\cO_p$ is greater than or equal to
the dimension of the fiber at $p$ of
$[(\pi_k)_*(\widetilde\cE_k)]'/(\pi_k)_*(\widetilde\cE_k)$.
But now we are in the situation of Lemma \ref{grassfiber}: this dimension
is at least the number of
negative line bundles in the direct sum decomposition of the pullback
of the universal subbundle of $OG$ under some positive-degree
map from a copy of $\bP^1_k$ to $OG$, and
this must be at least 2.
\end{proof}

\subsection{Degeneracy loci}
\label{degloci}
Degeneracy loci for vector bundles in type $D$
were defined using rank inequalities in \cite{KTdsp}.

\begin{defn}
\label{defnloci}
The degeneracy loci $W_\l$ and $W_\l(p)$ ($\l\in \D_n$, with $\ell=\ell(\l)$,
and $p\in \bP^1$)
are the following subschemes of $\bP^1\times OQ_d$:
\begin{gather*}
W_\l = \{\,x\in \bP^1\times OQ_d\,|\,
\rk(\E\to \cO\otimes V/F^{(i)\perp}_{n+1-\l_i})_x\lequ n+1-i-\l_i,
 i=1,\ldots,\ell+1 \,\}, \\
W_\l(p) = W_\l\cap (\{p\}\times OQ_d)
\end{gather*}
\end{defn}

Define also
$$h(n,d)=n(n+1)/2+2nd,$$
which is the dimension of the orthogonal Quot scheme $OQ_d$ when $d$
is a nonnegative integer.
As in types $A$ and $C$, we establish a Moving Lemma, and deduce from 
this that all
three-term Gromov--Witten invariants on $OG$ count points
in intersections of degeneracy loci on $OQ_d$.
\begin{movlemma}
Let $k$ be a positive integer, and let $p_1$, $\ldots$, $p_k$ be
distinct points on $\bP^1$.
Let $\l^1$, $\ldots$, $\l^k$ be partitions in $\D_n$, and
let us take the degeneracy loci $W_{\l^1}(p_1)$, $\ldots$, $W_{\l^k}(p_k)$
to be defined by isotropic flags of vector spaces in general position.
Consider the intersection
$$Z:=W_{\l^1}(p_1)\cap\cdots\cap W_{\l^k}(p_k).$$
Then $Z$ has dimension at most $h(n,d)-\sum_{i=1}^k |\l^i|$.
Moreover, $Z\cap OM_d$ is either empty or generically reduced and
of pure dimension $h(n,d)-\sum_i |\l^i|$;
also, $Z\cap (OQ_d\smallsetminus OM_d)$ has dimension at most
$h(n,d)-\sum_{i=1}^k |\l^i|-1$.
\end{movlemma}

The following are immediate consequences of the Moving Lemma.

\begin{cor}
\label{corgw}
Let $p$, $p'$, $p''\in \bP^1$ be distinct points.
Suppose $\l$, $\m$, $\n\in \D_n$ satisfy
$|\l|+|\m|+|\n|=h(n,d)$.
With degeneracy loci defined with respect to isotropic flags in general
position, the intersection $W_\l(p)\cap W_\m(p')\cap W_\n(p'')$
consists of finitely many reduced points, all contained in $OM_d$, and 
the corresponding Gromov--Witten invariant on $OG$ satisfies
$$\langle \t_\l, \t_\m, \t_\n \rangle_d = \#\bigl(W_\l(p)\cap W_\m(p')\cap
W_\n(p'')\bigr).$$
\end{cor}

\begin{cor}
\label{corzero}
If $p$ and $p'$ are distinct points of $\bP^1$ and if
$|\l|+|\m|=h(n,d)$, then
$W_\l(p)\cap W'_\m(p')=\emptyset$ for a general translate
$W'_\m(p')$ of $W_\m(p')$.
\end{cor}

The Moving Lemma itself is proved using an analysis of the boundary of
$OQ_d$. As in \cite{Ber} and \cite{KTlg}, this boundary is covered by
Grassmann bundles over smaller Quot schemes.

\begin{defn}
For $c\in (1/2)\Z$, with $c\gequ 1$,
we let $\pi_c\colon G_c\to \bP^1\times OQ_{d-c}$
denote the Grassmann bundle of $(2c)$-dimensional
quotients of the universal bundle $\E$ on $\bP^1\times OQ_{d-c}$.
The morphism $\beta_c\colon G_c\to OQ_d$ is given by the
modification of the sheaf sequence $\cE\to \cO\otimes V$
along the graph of the projection to $\bP^1$.
Precisely: let $\F_c$ denote the universal quotient bundle on $G_c$;
if $i_c$ denotes the morphism $G_c\to \bP^1\times G_c$ given by
$({\rm pr}_1\circ\pi_c,{\rm id})$, then $\E_c$ is defined as the kernel
of the natural morphism of sheaves
$({\rm id}\times ({\rm pr}_2\circ \pi_c))^*\E\to i_{c{*}}\pi_c^*\E$
composed with $i_{c{*}}$ applied to the morphism
to $\F_c$.
\end{defn}

\noindent
We also consider degeneracy loci with respect to the bundles $\E_c$.

\begin{defn}
We define $\WW_{c,\l}$ and $\WW_{c,\l}(p)$ to be
the following subschemes of $G_c$:
\begin{gather*}
\WW_{c,\l} = \{\,x\in G_c\,|\,
\rk(\E_c\to \cO\otimes V/F^\perp_{n+1-\l_i})_x\lequ n+1-i-\l_i,
i=1,\ldots,\ell+1\,\}, \\
\WW_{c,\l}(p) = \WW_{c,\l}(p) \cap \pi_c^{-1}(\{p\}\times OQ_{d-c})
\end{gather*}
\end{defn}

\subsection{Boundary structure of $OQ_d$}
\label{structure}
The boundary of $OQ_d$ is made up of points
where $\E\to\cO\otimes V$ drops rank at one
or more points of $\bP^1$;
note that wherever it drops rank, it does so by at least two
(by our definition of the Quot scheme).

\begin{thm}
\label{structurethm}
For any $d\in (1/2)\Z$, with $d\gequ 0$ and $d\ne 1/2$, we have
$$\dim OQ_d=\begin{cases} h(n,d) & \text{if $d\in \Z$}, \\
h(n,d)-5 & \text{otherwise}. \end{cases}$$
Furthermore, for
$c\in (1/2)\Z$, $c\gequ 1$,
the map $\beta_c\colon G_c\to OQ_d$ satisfies \\
{\em (i)} Given $x\in OQ_d$, if $\QQ_x$ has rank at least $n+1+c$ at
$p\in \bP^1$, then $x$ lies in the image of $\beta_c$. \\
{\em (ii)} The restriction of $\beta_c$ to
$\pi_c^{-1}(\bP^1\times OM_{d-c})$ is a locally closed immersion. \\
{\em (iii)} We have
$$\beta_c^{-1}(W_\l(p))=\pi_c^{-1}(\bP^1\times W_\l(p))\cup
\WW_{c,\l}(p)$$
where on the right, $W_\l(p)$ denotes the degeneracy locus in
$OQ_{d-c}$.
\end{thm}

The proof of Theorem \ref{structurethm}, as well as that of the Moving Lemma
(which uses Theorem \ref{structurethm}), is similar to that of 
the corresponding results in \cite{Ber} and \cite{KTlg}. 
Details are left to the reader.

\section{Intersection Theory on $OQ_d$}
\label{itoqd}

The Chow group of algebraic cycles modulo rational equivalence
of a scheme $\X$ is denoted $A_*\X$. We 
also employ the following notation.

\begin{defn}
Let $p$ denote a point of $\bP^1$. \\
(i) $\ev^p\colon OM_d\to OG$ is the evaluation at $p$ morphism; \\
(ii) $\tau(p)\colon OQ_d(p) \to OQ_d$ is the projection from the
{\em relative orthogonal Grassmannian}
$OQ_d(p) := OG_{n+1}(\QQ|_{\{p\}\times OQ_d})$, that is,
the closed subscheme
of the Grassmannian ${\rm Grass}_{n+1}$ of rank-$(n+1)$ quotients \cite{Groth}
of the indicated coherent sheaf,
defined by isotropicity and parity conditions on the kernel of the
composite morphism from $\cO_{{\rm Grass}_{n+1}}\otimes V$ to the universal
quotient bundle of the relative Grassmannian; \\
(iii) $\ev(p)\colon OQ_d(p)\to LG$ is the evaluation morphism
on the relative orthogonal Grassmannian; \\
(iv) $\ev_c^p\colon \pi_c^{-1}(\{p\}\times OM_{d-c})\to
OG(n+1-2c,2n+2)$ is evaluation at $p$.
\end{defn}

\begin{lemma}[\cite{KTlg}]
\label{simplepullbackprop}
Let $T$ be a projective variety which is a homogenous space
for an algebraic group $G$.
Let $\X$ be a scheme, equipped with an action of the group $G$.
Let $U$ be a $G$-invariant integral open subscheme of $\X$, and let
$f\colon U\to T$ be a $G$-equivariant morphism.
Then the map on algebraic cycles
\[
[V]\mapsto \bigl[\,f^{-1}(V)\overline{\phantom j}\,\bigr]
\]
respects rational equivalence, and hence induces a map on Chow groups
$A_*T\to A_*\X$.
\end{lemma}

\begin{cor}
\label{neededratequivcor}
Fix distinct points $p$, $p'\in \bP^1$.
For any $\l\in\D_n$ of length $\ell=\ell(\l)\gequ 3$,
the following cycles are rationally equivalent to zero on $OQ_d$ and
on $OQ_d(p')$:\\
{\em(i)}
$\bigl[\,(\ev^p)^{-1}(\X_\l)\overline{\phantom j}\,\bigr]-
\sum_{j=1}^{r-1} (-1)^{j-1}
\bigl[\,(\ev^p)^{-1}(\X_{\l_j,\l_r}\cap
\X'_{\l\smallsetminus\{\l_j,\l_r\}})\overline{\phantom j}\,\bigr].$ \\
{\em (ii)}
$\sum_{j=1}^{r-1} (-1)^{j-1}
\bigl[\,\beta_1\bigl((\ev_1^p)^{-1}(\Y_{\l_j,\l_r}\cap
\Y'_{\l\smallsetminus\{\l_j,\l_r\}})\bigr)\overline{\phantom j}\,\bigr].$ 
\end{cor}
\noindent
Here, and in the sequel, $\X'_{\m}$ and $\Y'_{\m}$ denote the translates
of $\X_{\m}$ and $\Y_{\m}$ by a general element of the group $SO_{2n+2}$.

\medskip

As is standard, for any closed subscheme $Z$ of a scheme $\X$,
$[Z]\in A_*\X$ denotes 
the class in the Chow group of the cycle associated to $Z$;
we let $[Z]_k$ be the dimension $k$ component of $[Z]$.

\begin{prop}
\label{pointmove}
{\em (a)} Suppose $\l$ and $\m$ are in $\D_n$, and let $p$, $p'$, $p''$ be
distinct points in $\bP^1$.
Assume that $\ell(\l)$ equals $1$ or $2$ and 
$\m$ has even length $\gequ 2$.
Let $k=h(n,d)-|\l|-|\m|$.
Then
\begin{align*}
\bigl[W_\l(p)\cap W'_\m(p')\bigr]_k
&= \bigl[W_\l(p)\cap W'_\m(p) \bigr]_k {\it \ in}\ A_*OQ_d,
                                                     \\
\bigl[\tau(p'')^{-1}\bigl(W_\l(p)\cap W'_\m(p')\bigr)\bigr]_k
&= \bigl[\tau(p'')^{-1}\bigl(W_\l(p)\cap W'_\m(p) \bigr)\bigr]_k
  {\it \ in}\ A_*OQ_d(p''),
\end{align*}
where $W'_{\m}(p)$
denotes degeneracy locus with respect to a general translate of
the isotropic flag of subspaces.

\medskip
\noindent
{\em (b)} In $A_* OQ_d$, we have
\begin{equation}
\label{twocycles}
\bigl[W_\l(p) \cap W'_\m(p)\bigr]_{k} 
= \bigl[\,(\ev^p)^{-1}(\X_\l\cap
\X'_\m)\overline{\phantom j}\,\bigr]
+ \bigl[\,\beta_1\bigl((\ev_1^p)^{-1}(\Y_\l\cap
\Y'_\m)\bigr)\overline{\phantom j}\,\bigr]
\end{equation}
and in $A_* OQ_d(p'')$, the cycle class
$\bigl[\tau(p'')^{-1}\bigl(W_\l(p)\cap W'_\m(p) \bigr)\bigr]_k$
is equal to the right-hand side of {\em (\ref{twocycles})}.
\end{prop}

\begin{proof}
By a dimension count which uses Proposition \ref{meetsmallergrass},
the irreducible components of dimension $k$
in $W_\l(p)\cap W'_\mu(p)$ are the ones indicated on the right-hand
side of (\ref{twocycles}).
As in \cite{KTlg}, now, the result follows from the rational equivalence
$\{p\}\sim\{p'\}$ on $\bP^1$, pulled back to
$Y:=(\bP^1\times W_{\l}(p))\cap W'_{\m}$
(or further pulled back to $OQ_d(p'')$), once we know that
the irreducible components of
$W_\l(p)\cap W'_\mu(p)$ of dimension $k$ are generically smooth
and in the closure of the complement of the fiber of $Y$ over $p$
(and that this remains true after pullback by $\tau(p'')$).
The `in the closure' portion of the claim follows by an argument involving
the Kontsevich compactification of $OM_d$, as in op.\ cit.
Generic smoothness is clear for
$(\ev^p)^{-1}(\X_\l\cap\X'_\m)$.
Transverality of a general translate also establishes generic smoothness
for the other component, once we notice that
any point $x$ in a dense open subset of
$\beta_1((\ev_1^p)^{-1}(\Y_\l\cap \Y'_\m))$
has the property that for any local $\C$-algebra $R$ with
residue field $R/{\mathfrak m}\simeq \C$ and any
$\psi\colon R\to W_\l(p)\cap W'_\m(p)$ with closed point mapping to $x$,
the map $\psi$ factors through the restriction of $\beta_1$ to
$\pi_1^{-1}(\{p\}\times OM_{d-1})$.

This assertion follows from elementary linear algebra,
but because of some tricky cases involving parity, we give a sketch of the
argument. Fix a basis $\{v_i\}$ of $V$ so that the symmetric form is
given by $\langle v_i, v_j\rangle = \delta_{i+j,2n+3}$. Without
loss of generality, the two general-position flags are
$$F_i = \Span(v_1, \ldots, v_i)$$
and
$$G^{(0)}_i = \Span(v_{2n+3-i}, \ldots, v_{2n+2}),$$
where the latter specifies $G_{n+1}$ or $\wt{G}_{n+1}$ equal to
$\Span(v_{n+2},\ldots,v_{2n+2})$ according to parity;
see (\ref{defnfi}).
We will show that the condition on $x$ holds whenever $x$ is in the
preimage of the intersection of the Schubert {\em cells} corresponding to
$\Y_\l$ and $\Y'_\m$, subject to the further condition that the
line on $OG$ parametrized by the point in $OG(n-1,2n+2)$ is incident to
$\X_\l$ and $\X'_\m$ at two {\em distinct} points.

Consider first the case $\ell(\l)=1$.
Let $x$ correspond to $(n-1)$-dimensional $A\subset V$ at the point $p$.
The condition to be in the Schubert cell for $\Y_\l$ implies that
$A\cap F_n^\perp=0$,
so $\rk(A\to V/F^{(i)}_{n+1})=n-1$ for any $i$.
By Definition \ref{defnloci},
the sheaf sequence corresponding to $\psi$ satisfies the rank condition
\begin{equation}
\label{rk1}
\rk(\cE\to \cO\otimes V/F^{(0)}_{n+1})\lequ n-1.
\end{equation}

Turning to the conditions coming from $\m$, we have
$\rk(A\cap G_{n+1}^{(1)})=n-\ell$, from membership in
the Schubert cell.
Suppose $n$ is even, so that $F^{(0)}_{n+1}=\wt{F}_{n+1}$ and
$G^{(1)}=G_{n+1}$ are disjoint.
Note that in this case Definition \ref{defnloci} imposes the condition
\begin{equation}
\label{rk2}
\rk(\cE\to \cO\otimes V/G_{n+1})\lequ n-\ell.
\end{equation}

The following basic argument is used to show that $\psi$ factors
through the restriction of $\beta_1$ to $\pi_1^{-1}(\{p\}\times OM_{d-1})$.
We have a sheaf sequence on $\bP^1_R$; after restricting to
$\AA^1_R$ the sheaf $\cE$ can be trivialized, so let us assume the
map to $\cO\otimes V$ is given by the $(2n+2)\times(n+1)$ 
matrix $L$ with values in $R[t]$,
with coordinates assigned so the top half of the matrix
corresponds to $\wt{F}_{n+1}$ and the bottom half corresponds to
$G_{n+1}$.
We may assume $t=0$ defines $p$, and also assume that mod ${\mathfrak m}$,
the rightmost two columns of $L$ vanish at $t=0$.
We localize at ${\mathfrak m}+tR[t]$.
It suffices to show that conditions (\ref{rk1}) and (\ref{rk2})
imply, after column operations, that the rightmost two columns of
$L$ have values in the ideal generated by $t$.
We have $\rk(A\to V/F_{n+1})=n-1$, that is,
some $(n-1)\times(n-1)$ minor in
the bottom half of $L$ has full rank. Now by performing column operations
and invoking (\ref{rk1}) we have all the entries in
the bottom right $(n+1)\times 2$ submatrix of $L$ lying in the ideal $(t)$.
Let $L'$ denote the top right $(n+1)\times 2$ submatrix of $L$.
The remaining isotropicity and rank conditions amount to $UL'=0$ mod $t$
for some matrix $U$, whose entries are polynomial functions of
the entries of $L$ in the first $n-1$ columns.
The condition that the line corresponding to $A$ meets the Schubert
varieties in distinct points implies that the nullspace of $U$ is trivial, 
and hence $L'$ has entries in $(t)$ as well.

If, instead, $n$ is odd, we use the fact that
$\rk(A\cap G_{n+1})=n+1-\ell$
(also a condition to be in the Schubert cell).  
{}From Definition \ref{defnloci},
\begin{equation}
\label{rk3}
\rk(\cE\to \cO\otimes V/G_{n+1})\lequ
\rk(\cE\to \cO\otimes V/G_n^\perp)\lequ n+1-\ell.
\end{equation}
Now $F^{(0)}_{n+1}=F_{n+1}$ and $G_{n+1}$ are disjoint, and the
basic argument applies, using (\ref{rk1}) and (\ref{rk3}).

In case $\ell(\l)=2$, we have $A\cap F_{n+1}^{(0)}=0$ and
(\ref{rk1}) still holds, so the argument is the same.
\end{proof}

\medskip
We now establish the rational equivalences on $OQ_d$ --- and on
$OQ_d(p'')$ --- which directly imply the quantum Giambelli formula
of Theorem \ref{ogthm}.

\begin{prop}
\label{quantumgiambellithm}
Fix $\l\in \D_n$ with $\ell=\ell(\l)\gequ 3$.
Set $r=2\lfloor(\ell+1)/2\rfloor$.
Let $p$, $p'$, $p''$ denote distinct points in $\bP^1$.
Then we have the following identity of cycle classes
\begin{equation}
\label{quantumgiambellieqone}
\bigl[\,(\ev^p)^{-1}(\X_\l)\overline{\phantom j}\,\bigr]=
\sum_{j=1}^{r-1} (-1)^{j-1}
\bigl[\,\bigl((\ev^p)^{-1}(\X_{\l_j,\l_r})\cap
(\ev^{p'})^{-1}(\X'_{\l\smallsetminus\{\l_j,\l_r\}})\bigr)
\overline{\phantom j}\,\bigr],
\end{equation}
both on $OQ_d$ and on $OQ_d(p'')$, where
$\X'_\m$ denotes the translate of $\X_\m$ by a generally chosen element of
the group $SO_{2n+2}$.
\end{prop}

\begin{proof}
Combining parts (a) and (b) of Proposition \ref{pointmove} gives
\begin{multline*}
\bigl[\,\bigl((\ev^p)^{-1} (\X_{\l_j,\l_r}) \cap
(\ev^{p'})^{-1}(\X'_{\l\smallsetminus\{\l_j,\l_r\}})\bigr)
\overline{\phantom j}\,\bigr] \\
= \bigl[\,(\ev^p)^{-1}(\X_{\l_j,\l_r}\cap
\X'_{\l\smallsetminus\{\l_j,\l_r\}})\overline{\phantom j}\,\bigr] +
\bigl[\,\beta_1\bigl((\ev_1^p)^{-1}(\Y_{\l_j,\l_r}\cap
\Y'_{\l\smallsetminus\{\l_j,\l_r\}})\bigr)\overline{\phantom j}\,\bigr]
\end{multline*}
for each $j$, with $1\lequ j\lequ r-1$.
Now (\ref{quantumgiambellieqone}) follows by summing and applying
(i) and (ii) of Corollary \ref{neededratequivcor}.
\end{proof}

\medskip

\begin{thm}
\label{qgthm}
Suppose $\l\in \D_n$, with $\ell=\ell(\l)\gequ 3$, and
set $r=2\lfloor(\ell+1)/2\rfloor$.
Then we have the following identity in $QH^*(OG)$:
\begin{equation}
\label{quantumgiambellieqtwo}
\tau_\l=\sum_{j=1}^{r-1} (-1)^{j-1} \tau_{\l_j,\l_r} 
\tau_{\l\smallsetminus\{\l_j,\l_r\}}.
\end{equation}
\end{thm}

\begin{proof} The classical component of (\ref{quantumgiambellieqtwo})
follows from the classical Giambelli formula for $OG$. To handle the 
remaining terms, apply a refined cap
product operation \cite[\S 8.1]{F}
along $\ev(p'')$ to general translates of $\X_\m$ for
all $\m\in\D_n$ with $|\m|=h(n,d)-|\l|$, and invoke
Corollaries \ref{corzero} and \ref{corgw} (as in the
proof of \cite[Thm.\ 5]{KTlg}).
\end{proof}

\section{Quantum Schubert calculus}
\label{qsc}

Our aim in this Section is to use Theorem \ref{ogthm} and the algebra of
$\wt{P}$-polynomials to find 
combinatorial rules that compute some of
the quantum structure constants that
appear in the quantum product of two Schubert classes.

\subsection{Algebraic background}
Let $\E_n$ denote the set of all partitions $\l$ with $\l_1\lequ n$.
The main properties of $\wt{Q}$-polynomials that we need are 
collected in \cite[\S 2.1 and \S 6.1]{KTlg}. They imply corresponding
facts about the $\wt{P}$-polynomials, in particular, that the
set $\{\wt{P}_{\l}(X)\ |\ \l\in\E_n\}$ is a free $\Z$-basis of the
ring $\Lambda_n'$ that they span. Hence, there exist 
integers $f(\l,\m;\,\n)$ such that
\begin{equation}
\label{Pmult}
\wt{P}_{\l}(X)\,\wt{P}_{\mu}(X)=\sum_{\n}f(\l,\m;\,\n)\,\wt{P}_{\n}(X);
\end{equation}
the constants $f(\l,\m;\,\n)$ are independent of $n$, and defined for
any $\l,\m,\nu\in\E_n$. The corresponding coefficients $e(\l,\m;\,\n)$
in the expansion of the product $\wt{Q}_{\l}(X)\,\wt{Q}_{\mu}(X)$
are related to these by the equation
\begin{equation}
\label{etof}
e(\l,\m;\,\n)=2^{\ell(\l)+\ell(\m)-\ell(\n)}f(\l,\m;\,\n).
\end{equation}

There are explicit combinatorial rules (involving signs in general)
for computing the integers $f(\l,\m;\,\n)$, which follow
from corresponding formulas for decomposing products of 
Hall-Littlewood polynomials; for more details, see \cite[\S 6.1]{KTlg}.
Define the connected components of a skew Young diagram
by specifying that two boxes are
connected if they share a vertex or an edge. We then have the
following Pieri type formula for $\l$ {\em strict}:
\begin{equation}
\label{Ppieri}
\wt{P}_{\l}(X)\,\wt{P}_k(X)=\sum_{\m} 2^{N'(\l,\m)}\,\wt{P}_{\m}(X),
\end{equation}
where the sum is over all partitions $\m\supset\l$ with
$|\m|=|\l|+k$ such that 
$\m/\l$ is a horizontal strip, and $N'(\l,\m)$ is one less than
the number of connected components of $\m/\l$. In particular, we
have $\wt{P}_{\l}(X)\wt{P}_n(X)=\wt{P}_{(n,\l)}(X)$ for all
$\l\in\D_n$.

When $\l$, $\mu$ and $\nu$ are strict partitions, the $f(\l,\mu;\,\n)$
are classical structure constants for $OG(n+1,2n+2)$,
\[
\t_{\l}\t_{\m}=\sum_{\n\in{\D_n}} f(\l,\mu;\,\n)\,\t_{\n},
\]
and hence are nonnegative integers. In this case,
Stembridge \cite{St} has given a combinatorial rule for the
numbers $f(\l,\mu;\,\n)$, analogous to the usual
Littlewood-Richardson rule in type $A$. Specifically,
$f(\l,\m;\,\n)$ is equal to the number of
marked tableaux of weight
$\l$ on the shifted skew shape ${\mathcal S}(\n/\m)$ satisfying certain
conditions (see \cite{St} and \cite[Sect.\ 6]{P} for more details).

\subsection{Quantum multiplication}
Recall from the Introduction that for any 
$\l,\m\in\D_n$ there is a formula 
\[
\t_{\l}\cdot \t_{\m} = 
\sum  f_{\l\m}^{\n}(n)\,\t_{\n}\, q^d
\]
in $QH^*(OG(n+1,2n+2))$, with each $f_{\l\m}^{\n}(n)$ equal to a
Gromov--Witten invariant
$\langle \t_{\l}, \t_{\m}, \t_{\wh{\n}} \rangle_d$ (defined when
$|\l|+|\m|=|\n|+2nd$). The nonnegative integer $f_{\l\m}^{\n}(n)$ 
counts the number of degree-$d$ rational maps
$\psi:\bP^1\ra OG$ such that $\psi(0)\in\X_{\l}$, $\psi(1)\in\X_{\m}$ and
$\psi(\infty)\in\X_{\wh{\n}}$, when the three Schubert varieties 
$\X_{\l}$, $\X_{\m}$ and $\X_{\wh{\n}}$ are in general position.

We adopt the convention that
$\t_{\l}=0$ for all non-strict partitions $\l$. 
Now Theorem \ref{ogthm} and the
Pieri rule (\ref{Ppieri}) give
\begin{cor}[Quantum Pieri Rule] 
\label{qupieri}
For any $\l\in\D_n$ and $k\gequ 0$ we have
\[
\t_{\l}\t_k=\sum_{\m} 2^{N'(\l,\m)}\t_{\m}+\sum_{\m\supset (n,n)}
2^{N'(\l,\m)} \t_{\m\ssm (n,n)} \, q
\]
where both sums are over $\mu\supset\l$ with
$|\mu|=|\l|+k$ such that $\mu/\l$ is a horizontal strip, 
and the second sum is restricted to those $\mu$ with two
parts equal to $n$.
\end{cor}
\noindent
In recent work with Buch \cite{BKT}, 
we give a more direct proof of the quantum Pieri rule for $OG$, and the
corresponding rule for the Lagrangian Grassmannian.

For any $d,n \gequ 0$ and partition $\n$, let $(n^d,\n)$ denote the
partition 
\[
(n,n,\ldots,n,\n_1,\n_2,\ldots),
\]
where $n$ appears $d$
times before the first 
component $\n_1$ of $\n$. Theorem \ref{ogthm} now gives

\begin{thm}
\label{qstructthm}
For any $d\gequ 0$ and strict partitions $\l,\m,\n\in\D_n$ 
with $|\n|=|\l|+|\m|-2nd$,
the quantum structure constant $f_{\l\m}^{\n}(n)$
satisfies $f_{\l\m}^{\n}(n)=f(\l,\m;\,(n^{2d},\n))$.
\end{thm}

We deduce that for any strict partitions $\l,\m,\n\in\D_n$, the
coefficient $f(\l,\m;\,(n^d,\n))$ is a nonnegative integer. The
constants $f(\l,\m;\,\n)$ can be negative; for example
\[
f(\rho_3,\rho_3;\,(4,4,2,2)) =-1.
\]
This follows from the Remark in \cite[\S 6.2]{KTlg}.

\subsection{The relation to $QH^*(LG(n-1,2n-2))$}
\label{LGconnection}
The quantum Pieri rule of Proposition \ref{qupieri} implies that
\[
\t_n\,\t_{\l}=
\begin{cases}
\t_{(n,\l)}&\text{if $\l_1<n$}, \\ \t_{\l\ssm (n)}\,q&\text{if $\l_1=n$}
\end{cases}
\]
in the quantum cohomology ring of $OG(n+1,2n+2)$.
Therefore, to compute all the Gromov--Witten invariants for
$OG$, it suffices to evaluate the
$\langle \t_{\l}, \t_{\m}, \t_{\n} \rangle_d$ for $\m,\n\in
\D_{n-1}$. Define a map $*:\D_n\ra\D_{n-1}$ by setting
$\l^*=(n-\l_{\ell},\ldots,n-\l_1)$ for
any partition $\l$ of length $\ell$, and $(0)^*=(0)$.

Partitions in $\D_{n-1}$ also parametrize the Schubert classes
$\s_{\l}$ in the (quantum) cohomology ring of the Lagrangian Grassmannian
$LG(n-1,2n-2)$, which was studied in \cite{KTlg}.
For the remainder of this paper, we let ${}':\D_{n-1}\ra\D_{n-1}$ denote
the duality involution for this space, so that the parts of $\l'$ 
complement the parts of $\l$ in the set $\{1,2,\ldots,n-1\}$. Notice 
that the restriction of $*$ to $\D_{n-1}$ defines a second involution
on this set, which was considered in \cite[\S 6.3]{KTlg}.

\begin{thm}
\label{OGLG}
Suppose that $\l\in\D_n$ is a non-zero 
partition with $\ell(\l)=2d+e+1$
for some nonnegative integers $d$ and $e$. For any
$\m,\n\in\D_{n-1}$, we have an equality
\begin{equation}
\label{lgogequ}
\langle \t_{\l}, \t_{\m}, \t_{\n} \rangle_d = 
\langle \s_{\l^*}, \s_{\m'}, \s_{\n'} \rangle_e
\end{equation}
of Gromov--Witten 
invariants for  $OG(n+1,2n+2)$
and $LG(n-1,2n-2)$, respectively. If $\l$ is zero or
$\ell(\l)<2d+1$, then
$\langle \t_{\l}, \t_{\m}, \t_{\n} \rangle_d =0$.
\end{thm}
\begin{proof}
Assume first that $\l_1<n$, so $\l\in\D_{n-1}$. We then have
\begin{align*}
\langle \t_{\l}, \t_{\m}, \t_{\n} \rangle_d &=
f(\l,\m;\, (n^{2d+1},\n')) \\
&= 2^{n+2d-\ell(\l)-\ell(\m)-\ell(\n)}\,e(\l,\m;\, (n^{2d+1},\n')) \\
&= 2^{n+4d+1-\ell(\l)-\ell(\m)-\ell(\n)}\,
\langle \s_{\l}, \s_{\m}, \s_{\n} \rangle_{2d+1}
\end{align*}
where the last equality comes from \cite[Thm.\ 6]{KTlg}.
The result now follows by applying the
eight-fold symmetry \cite[Thm.\ 7]{KTlg} for $QH^*(LG(n-1,2n-2))$,
which dictates
\begin{equation}
\label{lgsymmetry}
2^{n+2d} \,\langle \s_{\l}, \s_{\m}, \s_{\n} \rangle_{2d+1}
     =  2^{ \ell(\mu)+\ell(\nu)+e }\, \langle \s_{\l^*}, 
\s_{\mu'}, \s_{\nu'}\rangle_e.
\end{equation}
If $\l_1=n$, then
\[
\langle \t_{\l}, \t_{\m}, \t_{\n} \rangle_d=
\langle \t_{\l\ssm (n)}, \t_{\m}, \t_{(n,\n)} \rangle_d=
f(\l\ssm (n),\m;\, (n^{2d},\n')),
\]
and the previous analysis applies, since $\l^*=(\l\ssm (n))^*$.
\end{proof}

\medskip

Of course this theorem also provides an equality of Gromov--Witten
invariants going the other way. For any $\l,\m,\n\,\in\D_{n-1}$, we have
\[
\langle \s_{\l}, \s_{\m}, \s_{\n} \rangle_e = 
\begin{cases}
\langle \t_{\l^*}, \t_{\mu'}, \t_{\nu'}\rangle_d
&\text{if $\ell(\l)-e=2d+1$ is odd},\\
\langle \t_{(n,\l^*)}, \t_{\mu'}, \t_{\nu'}\rangle_d
&\text{if $\ell(\l)-e=2d$ is even}.
\end{cases}
\]
The $(\Z/2\Z)^3$-symmetry (\ref{lgsymmetry})
enjoyed by the Gromov--Witten invariants
for $LG(n-1,2n-2)$ implies a similar one for $QH^*(OG)$.

\begin{prop}
\label{ogsymmetry}
Let $\l\in\D_n$ be non-zero and $\m,\n\in\D_{n-1}$.
For any $d,e\gequ 0$ with $2d+e+1=\ell(\l)$, we have
\[
2^{\ell(\m)+\ell(\n)+e+\delta} 
\,\langle \t_{\l}, \t_{\m}, \t_{\n} \rangle_d
= 2^{n+2d}\,
\begin{cases}
\langle \t_{\l^*}, \t_{\mu'}, \t_{\nu'}\rangle_g
&\text{if $e=2g+1$ is odd},\\
\langle \t_{(n,\l^*)}, \t_{\mu'}, \t_{\nu'}\rangle_g
&\text{if $e=2g$ is even},
\end{cases}
\]
where $\delta=\delta_{\l_1,n}$ is the Kronecker symbol.
\end{prop}

We now obtain orthogonal analogues of 
\cite[Prop.\ 10]{KTlg} and \cite[Cor.\ 8]{KTlg}.

\begin{cor}
\label{whichnonzero}
Let $\l$, $\m$, $\n$ and $\delta$ be as in Proposition
\ref{ogsymmetry}.
Then the inequalities 
\begin{equation}
\label{ineq}
\ell(\m)+\ell(\n)-n+\delta \lequ 2 d \lequ \ell(\l)+\ell(\m)+\ell(\n)-n
\end{equation}
are necessary conditions for the Gromov--Witten invariant
$\langle \t_\l, \t_\m, \t_\n\rangle_d$ to be nonzero.
Moreover, if the two sides of either of the inequalities in
{\em(\ref{ineq})} differ by $0$ or $1$, then 
$\langle \t_\l, \t_\m, \t_\n\rangle_d$ is related by 
the eight-fold symmetry to a classical structure constant.
\end{cor}

\begin{cor}
\label{someproducts}
For any $\l\in\D_n$, we have
\[
\t_{\l}\cdot\t_{\rho_{n-1}}=
\begin{cases}
\t_{{\l^*}'}\,q^d
&\text{if $\ell(\l)=2d$ is even}, \\
\t_{(n,{\l^*}')}\,q^d
&\text{if $\ell(\l)=2d+1$ is odd}.
\end{cases}
\]
in $QH^*(OG)$.
In particular, 
\[
\t_{\rho_n}\cdot\t_{\rho_n}=
\begin{cases}
\t_n\,q^{n/2}
&\text{if $n$ is even}, \\
q^{(n+1)/2}
&\text{if $n$ is odd}.
\end{cases}
\]
\end{cor}

\section{Appendix: An identity in $\wt{P}$-polynomials}

We give a proof of the following identity, which is used to simplify
a formula for degeneracy loci in type $D$ \cite{KTdsp}.
The proof uses the algebraic formalism of \S \ref{pfids}.

\begin{prop}
Let $X=(x_1,\ldots,x_n)$ be an $n$-tuple of variables, and consider also
$\wt{X}=(-x_1,x_2,\ldots,x_n)$
and $X'=(x_2,\ldots,x_n)$.
Then, for any $\l\in\E_n$ of length $\ell\gequ 1$ we have
\begin{equation}
\label{ptident}
\sum_{i=1}^\ell (-1)^{i-1} \wt{P}_{\l\ssm\{\l_i\}}(X)e_{\l_i}(X') =
\wt{P}_\l(\wt{X})+(-1)^{\ell+1}\wt{P}_\l(X).
\end{equation}
\end{prop}

\begin{proof}
By homogeneity, (\ref{ptident}) is equivalent to the identity
\begin{equation}
\label{ptidentalt}
\sum_{i=1}^\ell (-1)^{i-1} \wt{Q}_{\l\ssm\{\l_i\}}(X)\wt{Q}_{\l_i}(X') =
\frac{1}{2}(\wt{Q}_\l(\wt{X})+(-1)^{\ell+1}\wt{Q}_\l(X)).
\end{equation}
To establish (\ref{ptidentalt}), we use identity (\ref{extnone}) and
are reduced to 
$$\sum_{i=1}^\ell (-1)^{i-1} 
\wt{Q}_{\l_i}(X')\sum_{\m\in B(\l\ssm\{\l_i\},k)}\wt{Q}_\m(X')=
\begin{cases}
\sum_{\m\in B(\l,k)}\wt{Q}_\m(X'), & \text{if $k\neq\ell$ mod $2$},\\
0 & \text{if $k=\ell$ mod $2$},
\end{cases}$$
for all integers $k$, where $B(\l,k)$ is defined as in the proof of Proposition
\ref{Boxprop}.
This corresponds to an identity in the algebra 
${\mathcal A}$ of formal variables with
imposed relations of \cite[\S 2.3]{KTlg}, which is similar to the algebra
${\mathcal B}$ of \S \ref{pfids}, except that only single bars appear.

Using the equalities
\begin{equation}
\label{newident}
[a,b](c)-[a,c](b)+[b,c](a)=0
\end{equation}
and
\begin{equation}
\label{newident2}
[a,b](\ov{c})-[a,c](\ov{b})+[b,c](\ov{a})=0
\end{equation}
in ${\mathcal A}$, one can verify, for each combination of 
parities of $k$ and $\ell$,
that the corresponding identity in ${\mathcal A}$ is true (one case,
that of $k$ odd, $\ell$ even, uses also the identity (\ref{pf})).
For example, when $k$ is even and $\ell$ is odd, we need to show that
\begin{equation}
\label{kevenloddnts}
\sum_{i=1}^\ell (-1)^{i-1}(\l_i)
\sum_{\m\in B(\l\ssm\{\l_i\},k)}
\sum\epsilon(\m,\n)(\nu_1,\nu_2)\cdots(\nu_{\ell-2},\nu_{\ell-1}) =
\sum_{\n\in B(\l,k)}(\n)
\end{equation}
where the innermost sum on the left is over all 
$(\ell-2)(\ell-4)\cdots(1)$ ways to
write the set of entries of $\m$ as a union of pairs
$\{\nu_1,\nu_2\}\cup\cdots\cup\{\nu_{\ell-2},\nu_{\ell-1}\}$.
Using (\ref{newident}), the sum of the terms on the left hand side
which contain a pair with exactly one bar vanishes.
The remaining terms are seen, using (\ref{newident}) and (\ref{newident2}),
to be equal to the Pfaffian expansion of the
right-hand side of (\ref{kevenloddnts}).
\end{proof}

\end{document}